%%%%%%%%%%%%%%%%%%%%%%%%%%%%%%%%%%%%%%%%%%%%%%%%%%%%%%%%%%%%%%%%%%%%%%%%%%%%%
%                                                                           %
%           WEIGHTED NORM INEQUALITIES FOR POLYNOMIAL EXPANSIONS            %
%               ASSOCIATED TO SOME MEASURES WITH MASS POINTS                %
%                                                                           %
%            J.J. GUADALUPE, M. PEREZ, F.J. RUIZ AND J.L VARONA             %
%                                                                           %
%                  ULTIMA MODIFICACION: 5 DE ABRIL DE 1995                  %
%             ESCRITO EN AMS-TEX 2.1 CON EL FORMATO AMSPPT.STY              %
%                                                                           %
%%%%%%%%%%%%%%%%%%%%%%%%%%%%%%%%%%%%%%%%%%%%%%%%%%%%%%%%%%%%%%%%%%%%%%%%%%%%%

\input amstex % THIS LINE IS NOT NECESSARY IF TYPESET BY AMS-TEX
\documentstyle{amsppt} % NOT NECESSARY IF TYPESET BY AMS-TEX WITH AMSPPT.STY
%\magnification=\magstep1

\define\N{{\Bbb N}} % <--- BLACBOARD BOLD N FOR THE SET OF NATURAL NUMBERS
\define\R{{\Bbb R}} % <--- BLACBOARD BOLD R FOR THE SET OF REAL NUMBERS

\define\lra{\longrightarrow} %<--- JUST AN ABREVIATED NAME
\define\BMO{\operatorname{BMO}} %<--- BMO SPACE
\define\overB{\,\overline{\!B}} %<--- OVERLINED B

\define\x{\relbar\joinrel}                             % | THIS IS AN EXTRA
\define\extralongarrow{\mapstochar\x\x\x\x\rightarrow} % | LONG ARROW 

\topmatter

\title Weighted norm inequalities for polynomial expansions associated to
      some measures with mass points\endtitle

\rightheadtext{Weighted norm inequalities for polynomial expansions}

\author Jos\'e J. Guadalupe, Mario P\'erez,	Francisco J. Ruiz and
     Juan L. Varona\endauthor

\leftheadtext{J. J. Guadalupe, M. P\'erez, F. J. Ruiz and J. L. Varona}

\address Jos\'e J. Guadalupe.
    Dpto. de Matem\'aticas y Computaci\'on.
    Universidad de La Rioja.
    26004 Logro\~no, Spain.\endaddress

\address Mario P\'erez.
    Dpto. de Matem\'aticas.
    Universidad de Zaragoza.
    50009 Zaragoza, Spain.\endaddress
    
\email mperez\@cc.unizar.es\endemail

\address Francisco J. Ruiz.
    Dpto. de Matem\'aticas.
    Universidad de Zaragoza.
    50009 Zaragoza, Spain.\endaddress

\address Juan L. Varona.
    Dpto. de Matem\'aticas y Computaci\'on.
    Universidad de La Rioja.
    26004 Logro\~no, Spain.\endaddress

\email jvarona\@siur.unirioja.es\endemail



\thanks Research supported by DGICYT under grant PB93-0228-C02-02 and by
        IER.\endthanks

\abstract Fourier series in orthogonal polynomials with respect to a measure
$\nu$ on $[-1,1]$ are studied when $\nu$ is a linear combination of a
generalized Jacobi weight and finitely many Dirac deltas in $[-1,1]$. We prove
some weighted norm inequalities for the partial sum operators $S_n$, their
maximal operator $S^*$ and the commutator $[M_b, S_n]$, where $M_b$ denotes
the operator of pointwise multiplication by $b \in \BMO$. We also prove some
norm inequalities for $S_n$ when $\nu$ is a sum of a Laguerre weight on $\R^+$
and a positive mass on $0$.\endabstract

\endtopmatter

\document

\head 0. Introduction.\endhead

Let $\nu$ be a positive Borel measure on $\R$ with infinitely many points of
increase and such that all the moments
$$
   \int_\R x^n d\nu \qquad (n = 0, 1, \dots)
$$
are finite. For each suitable function $f$, let $S_n f$ denote the $n$-th
partial sum of the Fourier expansion of $f$ with respect to the system of
orthogonal polynomials associated to $d\nu$.

The uniform boundedness of the operators $S_n \colon L^p(d\nu) \lra L^p(d\nu)$
($1 < p < \infty$) and some weighted versions $u S_n(v^{-1} \cdot) \colon
L^p(d\nu) \lra L^p(d\nu)$ have been characterized when $\nu$ is:
\roster
\item"a)" a Jacobi or generalized Jacobi weight on  $[-1,1]$ (see \cite{P},
\cite{M 1}, \cite{B});

\item"b)" a Laguerre weight on $\R^+$ or an Hermite weight on $\R$ (see
\cite{AW}, \cite{M 2}, \cite{M 3}).
\endroster

This uniform boundedness is equivalent, in rather general settings, to the
$L^p$ convergence of $S_n f$ to $f$.

Let us consider for simplicity the unweighted case. For a generalized Jacobi
weight, not only the uniform boundedness of the operators $S_n$ has been
studied, but also that of the maximal operator $S^*$ defined by
$$
   S^*f(x) = \sup_n |S_n f(x)|
$$
(see \cite{B}). For some orthogonal systems which include Jacobi polynomials
and Bessel functions, the maximal operator $S^*$ has been considered by Gilbert
(\cite{G}) by means of transference theorems.

Obviously, the boundedness of $S^*$ on $L^p(d\nu)$implies the uniform
boundedness of $S_n$ (and the $L^p(d\nu)$ convergence of $S_n f$ to $f$). But
it also implies, by standard arguments, the $\nu$-a.e. convergence of $S_n f$
to $f$. For these weights, the typical situation is that the operators $S_n
\colon L^p(d\nu) \lra L^p(d\nu)$ and even $S^*$ are uniformly bounded if and
only if $p$ belongs to some explicitly given open interval $(p_0, p_1)$ (the
interval of {\it mean convergence}). For short, in this case the $S_n$ are
said to be of {\it strong} $(p,p)$-type.

Then, for the endpoints $p = p_0, p_1$ of the interval of mean convergence it
is natural to study the {\it weak} $(p,p)$-type, i.e., the uniform boundedness
of the operators
$$
   S_n \colon L^p(d\nu) \lra L^{p,\infty}(d\nu),
$$
as well as the {\it restricted} {\it weak} $(p,p)$-type, i.e., the uniform
boundedness
$$
   S_n \colon L^{p,1}(d\nu) \lra L^{p,\infty}(d\nu).
$$
Here, $L^{p,r}(d\nu)$ stands for the classical Lorentz space of all measurable
functions $f$ satisfying
$$
   \eqalign{
   &\Vert f \Vert_{L^{p,r}(d\nu)} =
   \left(
   {r \over p} \int_0^\infty \left[ t^{1/p} f^*(t) \right]^r {dt \over t}
   \right)^{1/r} < \infty \quad (1 \leq p < \infty,\ 1 \leq r < \infty),
   \cr
   &\Vert f \Vert_{L^{p,\infty}(d\nu)} =
   \Vert f \Vert_{L^p_*(d\nu)} = \sup_{t > 0} t^{1/p} f^*(t) < \infty
   \quad (1 \leq p \leq \infty),
   \cr
   }
$$
where $f^*$ denotes the nonincreasing rearrangement of $f$. We refer the reader
to \cite{SW, Sect. V.3} for further information on these topics.

For $d\nu = dx$ on $[-1,1]$ (Fourier-Legendre series; $p_0 = 4/3$, $p_1 = 4$),
S. Chanillo (\cite{C}) proved that the $S_n$ are not of weak $(p,p)$-type for
$p = 4$ but they are of restricted weak $(p,p)$-type for $p = 4/3$ and $p =
4$. In \cite{GPV 1}, \cite{GPV 2} these results were established for any
Fourier-Jacobi series ($d\nu = (1-x)^\alpha (1+x)^\beta dx$ on $[-1,1]$,
$\alpha, \beta > -1$) and $p = p_0$, $p = p_1$. The $L^{p,r}$ behaviour of
$S_n$ was also studied by L. Colzani (\cite{Co}) for Fourier-Legendre series.

In this paper, we consider these problems for Fourier expansions with respect
to measures of the form
$$
   \nu = \mu + \sum_{i=1}^k M_i \delta_{a_i},
$$
where $M_i > 0$ ($i = 1, \dots, k$), $\delta_a$ denotes the Dirac delta on $a
\in \R$ and $\mu$ is a generalized Jacobi weight or, in some case, a Laguerre
or Hermite weight. In the particular case of a Jacobi weight and two mass
points on $1$ and $-1$, the corresponding orthonormal polynomials were studied
by Koornwinder in \cite{K} from the point of view of differential equations
(see also \cite{Ch}, \cite{AE}, \cite{Kr}, \cite{Li}). Our method consists of
relating the operators $S_n$ to some other operators similar to (and
expressible in terms of) the Fourier expansions with respect to $\mu$ and
polynomial modifications of $\mu$.

This method applies also to the commutator $[M_b, S_n]$, where $M_b$ is the
operator of pointwise multiplication by a given function $b$, i.e., $M_b f = b
f$. Given a linear operator $T$ acting on functions, say $T \colon L^p(d\nu)
\lra L^p(d\nu)$, and a function $b$, the commutator of $M_b$ and $T$ is
defined by
$$
   [M_b, T] f = b T(f) - T(b f).
$$
The first results on this commutator were obtained by R. R. Coifman, R.
Rochberg and G. Weiss (see \cite{CRW}). They proved that if $T$ is the
classical Hilbert transform and $1 < p < \infty$, then $[M_b, T]$ is a bounded
operator on $L^p(\R)$ if and only if $b \in \BMO(\R)$. The boundedness of this
commutator has been studied in more general settings by several authors
(see e.g. \cite{Bl}, \cite{ST 1}, \cite{ST 2}, \cite{ST 3}, \cite{GHST}).

Let us also mention that the boundedness of the commutator $[M_b, T]$ (or, in
our case, $[M_b, S_n]$) for $b$ in some real Banach space $B$ is closely
related to the analyticity (in our case, uniform analyticity) of the
operator-valued function
$$
   \eqalign{
   {\Cal T} \colon &\overB \lra {\Cal L}(L^p(d\nu), L^p(d\nu)) \cr
   &b \extralongarrow \ {\Cal T}(b) = M_{e^b} T M_{e^{-b}} \cr
   }
$$
in a neighbourhood of $0 \in \overB$, where $\overB$ denotes the
complexification of $B$ and ${\Cal L}(L^p(d\nu), L^p(d\nu))$ is the
space of bounded linear operators from $L^p(d\nu)$ into itself. In fact, the
first G\^ateaux differential of ${\Cal T}$ at $0$ in the direction $b \in B$,
is $$
   {d \over d s} {\Cal T}(s b) \Big|_{_{\scriptstyle s=0}} = [M_b, T]
$$
(see \cite{CM}, \cite{L} for further details). It is via this relationship that
R. R. Coifman and M. A. M. Murray proved (\cite{CM}) the uniform boundedness of
the commutator
$$
   [M_b, S_n] \colon L^2(d\nu) \lra L^2(d\nu)
$$
when $d\nu$ is a Jacobi weight ($d\nu = (1-x)^\alpha (1+x)^\beta dx$ on
$[-1,1]$) with $\alpha, \beta > -1/2$ and $b \in \BMO$.

Here, we prove the uniform boundedness of $[M_b, S_n]$ (as well as a weighted
version) in $L^p(d\nu)$, $1 < p < \infty$, where $d\nu$ is a generalized
Jacobi weight with possibly a finite collection of Dirac deltas on $[-1,1]$ and
again $b \in \BMO$.

This paper is organized as follows: in Sect. $1$ we present the basic notation
and technical results. In Sect. $2$ we consider the maximal operator $S^*$
related to a generalized Jacobi weight function with finitely many Dirac
masses on $[-1,1]$. As a consequence, the $L^p$ and a.e. convergence of the
Fourier series follow. For these measures (with some restriction), the
commutator $[M_b, S_n]$ is studied in Sect. $3$. Weak and restricted weak
boundedness at the endpoint of the interval of mean convergence for Jacobi
weights with Dirac masses on $[-1,1]$ are the subject of Sect. $4$. Finally, in
Sect. $5$ we point out how $L^p$ boundedness of Fourier expansions with respect
to a Laguerre or Hermite weight with a positive mass on $0$ can be established.

\head 1. Notations and technical results.\endhead

Let $\mu$ be a positive Borel measure on $\R$ with infinitely many
points of increase and such that all the moments
$$
   \int_\R x^n d\mu \qquad (n = 0, 1, \dots)
$$
are finite. Let $a_i \in \R$ ($i = 1, \dots, k$) with $a_i \not= a_j$ for $i
\not= j$ and assume $\mu(\{a_i\}) = 0$ ($i = 1, \dots, k$). Let $M_i > 0$ ($i
= 1, \dots, k$) and write
$$
   \nu = \mu + \sum_{i=1}^k M_i \delta_{a_i},
   \tag 1
$$
where $\delta_a$ denotes a Dirac delta on $a$:
$$
   \int_\R f d\delta_a = f(a).
$$
Then, there is a sequence $\{P_n\}_{n \geq 0}$ of polynomials,
$$
   P_n(x) = k_n x^n + \dots, \quad k_n > 0, \quad \deg P_n = n
$$
such that
$$
   \int_{\R} P_n P_m d\nu =
   \cases 0, &\text{if $n \not= m$;}\cr
          1, &\text{if $n=m$.}
   \endcases
$$
The $n$-th partial sum operator of the Fourier expansion in terms of $P_n$ is
the operator $S_n$ given by
$$
   S_n f(x) = \int_\R L_n(x,y) f(y) d\nu(y),
$$
where
$$
   L_n(x,y) = \sum_{j=0}^n P_j(x) P_j(y)
$$
is the $n$-th kernel relative to the measure $d\nu$. If we denote
$$
   T_n f(x) = \int_\R L_n(x,y) f(y) d\mu(y),
$$
then, according to (1), we have
$$
   S_n f(x) = T_n f(x) + \sum_{i=1}^k M_i L_n(x,a_i) f(a_i).
   \tag 2
$$
By a weight function we mean a non-negative, measurable function. We are
interested in finding conditions for the uniform boundedness of the operators
$$
   u S_n(v^{-1} \cdot) \colon L^p(d\nu) \lra L^p(d\nu)
$$
where $u$ and $v$ are weights, i.e., for the inequality
$$
   \Vert u S_n(v^{-1} f) \Vert_{L^p(d\nu)} \leq C \Vert f \Vert_{L^p(d\nu)}
$$
to hold for $n \geq 0$ and $f \in L^p(d\nu)$, and also for the weaker
boundedness
$$
   u S_n(v^{-1} \cdot) \colon L^p(d\nu) \lra L^{p,\infty}(d\nu)
$$
or
$$
   u S_n(v^{-1} \cdot) \colon L^{p,1}(d\nu) \lra L^{p,\infty}(d\nu).
$$
Actually, the last one is equivalent (see \cite{SW, Theorem 3.13}) to
$$
   \Vert u S_n(v^{-1} \chi_E) \Vert_{L^{p,\infty}(d\nu)}
   \leq C \Vert \chi_E \Vert_{L^p(d\nu)}
$$
for every measurable set $E$. In this context, notice that the values
$u(a_i)$, $v(a_i)$ are significant here, since $\nu(\{a_i\}) > 0$.

In what follows, given $1 \leq p \leq \infty$ we will denote by $p'$ the
conjugate exponent, i.e., $1 \leq p' \leq \infty$, $1/p + 1/p' = 1$. Also, we
will take $0 \cdot \infty = 0$ and by $C$ we will mean a constant, not
depending on $n$, $f$, but possibly different at each occurrence.

Then, we have the following results:

\proclaim{Lemma 1}
With the above notation, let $1 < p < \infty$, $1 <
q < \infty$, $1 \leq r \leq \infty$, $1 \leq s \leq \infty$; let $u,v$ be two
weight functions on $\R$ with $u(a_i) < \infty$, $0 < v(a_i)$, $i=1,\dots, k$.
Then, there exists some constant $C>0$ such that
$$
   \Vert u S_n(v^{-1}f) \Vert_{L^{p,r}(d\nu)} \leq
   C \Vert f \Vert_{L^{q,s}(d\nu)}
   \tag 3
$$
for every $f \in L^{q,s}(d\nu)$, $n \geq 0$ if and only if there exists $C>0$
such that:
\roster
\item"a)" $\Vert u T_n(v^{-1}f) \Vert_{L^{p,r}(d\mu)} \leq C \Vert f
\Vert_{L^{q,s}(d\mu)}, \qquad f \in L^{q,s}(d\mu), \ n \geq 0$.

\item"b)" $u(a_i) \Vert v^{-1} L_n(x,a_i) \Vert_{L^{q',s'}(d\mu)} \leq C,
\qquad n \geq 0,\ i = 1, \dots, k$.

\item"c)" $v(a_i)^{-1} \Vert u L_n(x,a_i) \Vert_{L^{p,r}(d\mu)} \leq C,
\qquad n \geq 0,\ i = 1, \dots, k$.
\endroster
\endproclaim

We can state a similar result about the maximal operator $S^*$ defined by
$$
   S^*f(x) = \sup_n |S_nf(x)|.
$$
Let us also take
$$
   T^*f(x) = \sup_n |T_nf(x)|
$$
and
$$
   L^*(x,y) = \sup_n |L_n(x,y)|.
$$

\proclaim{Lemma 2}
With the above notation, let $1 < p < \infty$, $1 <
q < \infty$, $1 \leq r \leq \infty$, $1 \leq s \leq \infty$; let $u,v$ be two
weight functions on $\R$ with $u(a_i) < \infty$, $0 < v(a_i)$, $i=1,\dots, k$.
Then, there exists some constant $C>0$ such that
$$
   \Vert u S^*(v^{-1}f) \Vert_{L^{p,r}(d\nu)} \leq
   C \Vert f \Vert_{L^{q,s}(d\nu)},
$$
for every $f \in L^{q,s}(d\nu)$ if and only if there exists $C>0$ such that:
\roster
\item"a)" $\Vert u T^*(v^{-1}f) \Vert_{L^{p,r}(d\mu)} \leq C \Vert f
\Vert_{L^{q,s}(d\mu)}, \qquad f \in L^{q,s}(d\mu)$.

\item"b)"$u(a_i) \Vert v^{-1} L_n(x,a_i) \Vert_{L^{q',s'}(d\mu)} \leq C,
\qquad n \geq 0,\ i = 1, \dots, k$.

\item"c)" $v(a_i)^{-1} \Vert u L^*(x,a_i) \Vert_{L^{p,r}(d\mu)} \leq C,
\qquad i = 1, \dots, k$.
\endroster
\endproclaim

Finally, we also have the analogous result for the commutator (notice that
$b(a_i)$ is significant here, too):

\proclaim{Lemma 3}
With the above notation, let $1 < p < \infty$, $1 < q <
\infty$, $1 \leq r \leq \infty$, $1 \leq s \leq \infty$; let $u,v$ be two
weight functions on $\R$ with $u(a_i) < \infty$, $0 < v(a_i)$, $i=1,\dots, k$.
Let $b$ be a function on $\R$ with $b(a_i) < \infty$, $i = 1, \dots, k$.
Then, there exists some constant $C>0$ such that
$$
   \Vert u [M_b, S_n](v^{-1}f) \Vert_{L^{p,r}(d\nu)} \leq
   C \Vert f \Vert_{L^{q,s}(d\nu)}
$$
for every $f \in L^{q,s}(d\nu)$, $n \geq 0$ if and only if there exists $C>0$
such that:
\roster
\item"a)" $\Vert u [M_b, T_n](v^{-1}f) \Vert_{L^{p,r}(d\mu)} \leq C \Vert f
\Vert_{L^{q,s}(d\mu)}, \qquad f \in L^{q,s}(d\mu), \ n \geq 0$.

\item"b)" $u(a_i) \Vert v^{-1} L_n(x,a_i) [b(x) - b(a_i)]
\Vert_{L^{q',s'}(d\mu)} \leq C, \qquad n \geq 0,\ i = 1, \dots, k$.

\item"c)" $v(a_i)^{-1} \Vert u L_n(x,a_i) [b(x) - b(a_i)] \Vert_{L^{p,r}(d\mu)}
\leq C, \qquad n \geq 0,\ i = 1, \dots, k$.
\endroster
\endproclaim

\remark{Remark}
From the definition, we have $\Vert \chi_E
\Vert_{L^{p,r}(d\sigma)} = \sigma(E)^{1/p}$ for any measure $\sigma$ and any
measurable set $E$. In particular, any function is a.e. a characteristic
function with respect to a measure of the form $M \delta_a$, thus $\Vert f
\Vert_{L^{p,r}(M\delta_a)} = M^{1/p} |f(a)|$. As a consequence, we obtain in
our case
$$
   \Vert f \Vert_{L^{p,r}(d\nu)} \sim \Vert f \Vert_{L^{p,r}(d\mu)} +
   \sum_{i=1}^k M_i^{1/p} |f(a_i)|,
$$
where ``$\sim$'' means that the ratio is bounded above and below by two
positive constants. Actually, the $M_i$ can be removed.
\endremark

\demo{Proof of Lemma 1}
I) Suppose (3) holds. Let
$f \in {L^{q,s}(d\mu)}$ and put $g(a_i) = 0$ ($i=1, \dots, k$), $g(x) = f(x)$
otherwise. Clearly $f = g$ ($\mu$-a.e.), so that from relation (2) and the
definition of $T_n$ it follows
$$
   u S_n(v^{-1}g) = u T_n(v^{-1}g) = u T_n(v^{-1}f).
$$
It is also easy to deduce
$$
   \Vert g \Vert_{L^{q,s}(d\nu)} = \Vert g \Vert_{L^{q,s}(d\mu)} =
   \Vert f \Vert_{L^{q,s}(d\mu)}
$$
($g$ has the same distribution function with respect to $\nu$ and $\mu$ and it
coincides with the distribution function of $f$ with respect to $\mu$).
Therefore, from (3) applied to $g$, we have
$$
   \Vert u T_n(v^{-1} f) \Vert_{L^{p,r}(d\nu)}
   \leq C \Vert f \Vert_{L^{q,s}(d\mu)}
   \tag 4
$$
for every $f \in {L^{q,s}(d\mu)}$ and $n \geq 0$.

Now, set $f = \chi_{\{a_i\}}$; then,
$$
   u(x) S_n(v^{-1}f)(x) = u(x) M_i L_n(x,a_i) v(a_i)^{-1},
$$
$$
   \Vert f \Vert_{L^{q,s}(d\nu)} = M_i^{1/q}
$$
(for the last equality, see the previous remark). Therefore, from (3),
applied to $f$, it follows
$$
   v(a_i)^{-1} \Vert u L_n(x,a_i) \Vert_{L^{p,r}(d\nu)} \leq C
   \tag 5
$$
for every $n \geq 0$ and $i = 1, \dots, k$.

II) Conversely, suppose (4) and (5) hold. Then, for every $f \in L^{q,s}(d\nu)$
we have, from (2),
$$
   \Vert u S_n(v^{-1}f) \Vert_{L^{p,r}(d\nu)}
$$
$$
   \leq
   \Vert u T_n(v^{-1}f) \Vert_{L^{p,r}(d\nu)} +
   \sum_{i=1}^k M_i \Vert u L_n(x,a_i) v(a_i)^{-1} f(a_i) \Vert_{L^{p,r}(d\nu)}
$$
$$
   \leq C \Vert f \Vert_{L^{q,s}(d\mu)} + C \sum_{i=1}^k M_i |f(a_i)|
   \leq C \Vert f \Vert_{L^{q,s}(d\nu)},
$$
using the previous remark. That is: (3) is equivalent to (4) and (5). We will
see now that (5) is the same as $c)$ and that (4) is equivalent to $a)$ and
$b)$.

III) Fix an $i \in \{1, \dots, k\}$. Then
$$
   \Vert u L_n(x,a_i) \Vert_{L^{p,r}(d\nu)}
$$
$$
   \sim \Vert u L_n(x,a_i) \Vert_{L^{p,r}(d\mu)} +
   \sum_{j=1}^k M_j^{1/p} u(a_j) |L_n(a_j,a_i)|.
$$
Now, from Cauchy-Schwarz inequality we have
$$
   |L_n(a_j,a_i)| \leq L_n(a_j,a_j)^{1/2} L_n(a_i,a_i)^{1/2}
$$
and the sequences $\{L_n(a_j,a_j)\}_{n \geq 0}$ ($j = 1, \dots, k$) are
bounded, since $\nu(\{a_j\}) > 0$ (see \cite{N, p. 4}). Hence,
$$
   \Vert u L_n(x,a_i) \Vert_{L^{p,r}(d\mu)}
   \leq \Vert u L_n(x,a_i) \Vert_{L^{p,r}(d\nu)} \leq
   \Vert u L_n(x,a_i) \Vert_{L^{p,r}(d\mu)} + C.
$$
This means that (5) is actually equivalent to $c)$.

IV) Let us take now condition (4). From the previous remark again,
$$
   \Vert u T_n(v^{-1}f) \Vert_{L^{p,r}(d\nu)}
$$
$$
   \sim \Vert u T_n(v^{-1}f) \Vert_{L^{p,r}(d\mu)} +
   \sum_{i=1}^k M_i^{1/p} u(a_i) |T_n(v^{-1}f)(a_i)|.
$$
Thus, (4) is equivalent to condition $a)$, together with:
$$
   u(a_i) |T_n(v^{-1}f)(a_i)| \leq C \Vert f \Vert_{L^{q,s}(d\mu)}.
$$
Having in mind that
$$
   u(a_i) T_n(v^{-1}f)(a_i)
   = u(a_i) \int_\R v(x)^{-1} L_n(x,a_i) f(x) d\mu(x),
$$
that means, by duality,
$$
   u(a_i) \Vert v(x)^{-1} L_n(x,a_i) \Vert_{L^{q',s'}(d\mu)} \leq C,
$$
i.e., condition $b)$. \qed
\enddemo

The proofs of Lemmas 2 and 3 are essentially the same, so we omit them.

The following result provides sufficient conditions for the uniform
boundedness of the operators $T_n$ in terms of the boundedness of the Fourier
series corresponding to the measure $\mu$ and other related measures.
Recalling that
$$
   \nu = \mu + \sum_{i=1}^k M_i\delta_{a_i},
$$
we define, for each set $A \subseteq \{a_1,\dots,a_k\}$, the measure
$$
   d\mu^A(x) = \prod_{a_i\in A}(x-a_i)^2 \ d\mu(x)
$$
(for $A=\emptyset$ we just get $d\mu$) and the associated partial sum
operators ${\widetilde S}_n^A$. We also define the weight
$$
   w^A(x) = \prod_{a_i\in A}|x-a_i|^{1-2/p}.
$$
With this notation, we have:

\proclaim{Lemma 4}
If for each $A \subseteq \{a_1,\dots,a_k\}$ there
exists a constant $C$ such that
$$
   \Vert u w^A {\widetilde S}_n^A([v w^A]^{-1} f) \Vert _{L^p(d\mu^A)}
   \leq C \Vert f \Vert _{L^p(d\mu^A)}
$$
for every $n \geq 0$ and $f \in L^p(d\mu^A)$, then there also
exists a constant $C$ such that
$$
   \Vert uT_n(v^{-1} f) \Vert _{L^p(d\mu)}
   \leq C \Vert f \Vert _{L^p(d\mu)}
$$
for $n \geq 0$ and $f \in L^p(d\mu)$.
\endproclaim

\demo{Proof}
Let us denote by $K_n^A(x,y)$ the $n$-th kernel relative to the
measure $d\mu^A$. In the case $k=1$, we have
$$
   L_n(x,y) = C_n K_n^\emptyset(x,y)
   + (1-C_n)(x-a_1)(y-a_1)K_{n-1}^{\{a_1\}}(x,y)
$$
with $0< C_n <1$, $\forall n \in \N$ (see \cite{GPRV 3, Proposition 5}). By
induction on $k$, it can be shown that
$$
   L_n(x,y) = \sum_A \ C_n^A \left[\ \prod_{a_i \in A} (x-a_i)(y-a_i)\right]
   \ K_{n-|A|}^A(x,y),
$$
where the sum is taken over all the subsets $A \subseteq \{a_1,\dots,a_k\}$,
$|A|$ is the cardinal of $A$ and for each $n$
$$
   \sum_A C_n^A = 1, \qquad 0 < C_n^A < 1 \quad \forall A.
$$
From this expression, we deduce
$$
   T_n(v^{-1} f)(x) = \sum_A \ C_n^A \left[\ \prod_{a_i \in A} (x-a_i)\right]
   \ {\widetilde S}_{n-|A|}^A
   \left({v(y)^{-1}f(y) \over \prod_{a_i \in A} (y-a_i)},x\right).
$$
Thus, for the uniform boundedness of the operators $T_n$ it is enough to have
$$
   \left\Vert u(x) \prod_{a_i \in A} (x-a_i)
   {\widetilde S}_n^A\left({v(y)^{-1}f(y) \over
   \prod_{a_i \in A}(y-a_i)},x\right)
   \right\Vert _{L^p(d\mu)}
   \leq C \Vert f \Vert_{L^p(d\mu)}
$$
for $n \geq 0$ and $f \in L^p(d\mu)$. This is simply our hypothesis, except
for a change of notation. \qed
\enddemo

Analogous results can be stated about the maximal operator $T^*$ and
the commutator.

\head 2. The maximal operator $S^*$ for a generalized  Jacobi weight with
mass points on $[-1,1]$.\endhead

Let $d\mu = w dx$, with $w$ be a generalized Jacobi weight, that is:
$$
   w(x) = h(x) (1-x)^\alpha (1+x)^\beta
   \prod_{i=1}^N |x-t_i|^{\gamma_i},\quad x \in [-1,1]
$$
where:
\roster
\item"a)"	$\alpha$,$\beta$,$\gamma_i>-1$, $t_i \in (-1,1)$, $t_i \not= t_j$ for
$i \not= j$;

\item"b)"	$h$ is a positive, continuous function on $[-1,1]$ and
$w(h,\delta)\delta^{-1} \in L^1(0,2)$, $w(h,\delta)$ being the modulus of
continuity of $h$.
\endroster

Let $\displaystyle \nu = \mu + \sum_{i=1}^k M_i \delta_{a_i}$, with $M_i > 0$,
$a_i \in [-1,1]$ ($i = 1, \dots, k$).

Let us take also two weights $u$ and $v$ defined on $[-1,1]$ as follows:
$$
   u(x) = (1-x)^a(1+x)^b \prod_{i=1}^N |x-t_i|^{g_i},
   \quad \text{ if } x \not= a_i \ \forall i; \quad 0 < u(a_i) < \infty;
$$
$$
   v(x) = (1-x)^A(1+x)^B \prod_{i=1}^N |x-t_i|^{G_i},
   \quad \text{ if } x \not= a_i \ \forall i; \quad 0 < v(a_i) < \infty,
   \tag 6
$$
where $a,b,g_i,A,B,G_i \in \R$.

\proclaim{Theorem 5}
Let $1<p<\infty$. Then, there exists a constant $C>0$ such that
$$
   \Vert u S^*(v^{-1}f) \Vert _{L^p(d\nu)}	\leq C \Vert f \Vert _{L^p(d\nu)}
$$
for every $f \in L^p(d\nu)$ if and only if the inequalities
$$
   \cases
   A + (\alpha +1) ({1 \over p} - {1 \over 2})
   < \min \{{1 \over 4},{\alpha +1 \over 2}\}
   \cr
   B + (\beta +1) ({1 \over p} - {1 \over 2})
   < \min \{{1 \over 4},{\beta +1 \over 2}\}
   \cr
   G_i + (\gamma_i +1) ({1 \over p} - {1 \over 2})
   < \min \{{1 \over 2},{\gamma_i +1 \over 2}\} \quad (i=1,\dots,N)
   \endcases
   \tag 7
$$
$$
   \cases
   a + (\alpha +1) ({1 \over p} - {1 \over 2})
   > - \min \{{1 \over 4},{\alpha +1 \over 2}\}
   \cr
   b + (\beta +1) ({1 \over p} - {1 \over 2})
   > - \min \{{1 \over 4},{\beta +1 \over 2}\}
   \cr
   g_i + (\gamma_i +1) ({1 \over p} - {1 \over 2})
   > - \min \{{1 \over 2},{\gamma_i +1 \over 2}\} \quad (i=1,\dots,N)
   \endcases
   \tag 8
$$
and
$$
   \vbox{\halign{\hfil #  & \qquad \hfil #  &  \qquad \hfil # \hfil \cr
   $A \leq a$, & $B \leq b$, & $G_i \leq g_i \quad (i=1,\dots,N)$ \cr
   }}
   \tag 9
$$
hold.
\endproclaim

\remark{Remark}
This result is true in the case of a generalized Jacobi
weight, with no mass points. It was proved by V. M. Badkov (see \cite{B}) for
one weight ($u = v$). In the two weight case, the ``if'' part can be obtained
as a consequence, by inserting a suitable weight $\rho$, $u \leq \rho \leq v$.
Regarding the ``only if'' part,
$$
   \Vert u S^*(v^{-1} f) \Vert _{L^p(d\nu)}	\leq C \Vert f \Vert _{L^p(d\nu)}
$$
implies
$$
   \Vert u S_n(v^{-1} f) \Vert _{L^p(d\nu)}	\leq C \Vert f \Vert _{L^p(d\nu)}
$$
and this, in turn, implies (7), (8) and (9) (see \cite{GPRV 1}).
\endremark

\demo{Proof of the theorem}
Assume
$$
   \Vert uS^*(v^{-1}f) \Vert _{L^p(d\nu)}	\leq C \Vert f \Vert _{L^p(d\nu)}
$$
for $f \in L^p(d\nu)$. Proceeding as in \cite{GPRV 1, Theorem 6}, we obtain
(7), (8) and (9).

Assume now that (7), (8) and (9) hold. According to Lemma 2 with $p = q = r =
s$ and the analog of Lemma 4, we only need to prove the following inequalities:
\roster
\item"a)" $\Vert u w^A ({\widetilde S}^A)^*([v w^A]^{-1} f)
   \Vert_{L^p(d\mu^A)}
   \leq C \Vert f \Vert _{L^p(d\mu^A)},
   \qquad f \in L^p(d\mu^A);$

\item"b)" $\Vert v^{-1}L_n(x,a_i) \Vert _{L^{p'}(w)} \leq C, \qquad n \geq 0,
\ i=1,\dots,k$;

\item"c)" $\Vert uL^*(x,a_i) \Vert _{L^p(w)} \leq C, \qquad i=1,\dots,k$.
\endroster

Condition $a)$ refers to the boundedness of the maximal operators related to
the measures $d\mu^A$, which are generalized Jacobi weights, with no mass
points. In this case, the corresponding inequalities $(7)$, $(8)$ and $(9)$,
with the appropriate exponents, imply the boundedness. It is easy to see that
they actually hold.

To check inequalities $b)$ and $c)$, we can use the following estimates for the
kernels $L_n(x,a_i)$ (see \cite{GPRV 3}): if $a_i \not= \pm 1$,
$$
   |L_n(x,a_i)| \leq C
   (1-x+n^{-2})^{-(2\alpha+1)/4} (1+x+n^{-2})^{-(2\beta+1)/4} 
   \prod_{t_j \not= a_i} (|x-t_j|+n^{-1})^{-\gamma_j/2};
   \tag 10
$$
if $1$ is a mass point,
$$
   |L_n(x,1)| \leq C (1+x+n^{-2})^{-(2\beta+1)/4} 
   \prod_{i=1}^N (|x-t_i|+n^{-1})^{-\gamma_i/2};
   \tag 11
$$
if $-1$ is a mass point,
$$
   |L_n(x,-1)| \leq C (1-x+n^{-2})^{-(2\alpha+1)/4}
   \prod_{i=1}^N (|x-t_i|+n^{-1})^{-\gamma_i/2}.
   \tag 12
$$
It is not difficult to see that these inequalities, together with (7), (8)
and (9), lead to $b)$ and $c)$. \qed
\enddemo

\proclaim{Corollary 6}
With the notation of Theorem 5, the uniform boundedness
$$
   \Vert uS_n(v^{-1} f) \Vert _{L^p(d\nu)}	\leq C \Vert f \Vert _{L^p(d\nu)}
$$
holds for every $n \geq 0$ and $f \in L^p(d\nu)$ if and only if the
inequalities (7), (8) and (9) are verified.
\endproclaim

\demo{Proof}
The ``if'' part follows directly from the theorem and the
``only if'' part can be proved again as in \cite{GPRV 1}. \qed
\enddemo

\proclaim{Corollary 7}
Let $v$ be a weight verifying (6). Then,
$$
   v(x) S_n (v^{-1}f) (x) \lra f(x), \qquad \nu\hbox{-a.e.}
$$
for every $f \in L^p(d\nu)$ if and only if the inequality (7) holds.
\endproclaim

\demo{Proof}
The ``only if'' part is a consequence of \cite{GPRV 2, Theorem 3}. For the
``if'' part, we can take a weight $u$ such that the pair $(u,v)$ satisfies the
conditions of Theorem $5$ and it follows by standard arguments. Notice that
the weight $u$ does not play any role for the almost everywhere convergence.
\qed
\enddemo

\head 3. The commutator $[M_b, S_n]$ for a generalized Jacobi weight with
mass points on $[-1,1]$.\endhead

In this section we will adopt the notation of Sect. $2$, with the additional
restriction $\gamma_i \geq 0$, $i = 1, \dots, N$. We will write $I = [-1,1]$.

The space $\BMO(I)$ (in the sequel, $\BMO$) consists of the functions (modulus
the constants) of bounded mean oscillation, i.e., with
$$
   \Vert b \Vert_{\BMO} = \sup_J {1 \over |J|} \int_J |b - b_J| < \infty,
$$
where the supremum is taken over all the intervals $J \subseteq I$, $|J|$ means
the Lebesgue measure of $J$ and
$$
   b_J = {1 \over |J|} \int_J b
$$
(the integrals are taken with respect to the Lebesgue measure). Given $1 < p
< \infty$, we also have (see \cite{GR})
$$
   \Vert b \Vert_{\BMO} \sim
   \sup_J \left({1 \over |J|} \int_J |b - b_J|^p \right)^{1/p}.
   \tag 13
$$
Given $1 < p < \infty$ and a weight $\phi$ in Muckenhoupt's class $A_p$, the
commutator $[H, M_b]$ of the Hilbert transform on $I$ is bounded in
$L^p(\phi)$ if and only if $b \in \BMO$ (see \cite{Bl}). The norm of the
commutator depends on the $A_p$ constant of $\phi$. We refer the reader to
\cite{GR} for further references on $A_p$ weights.

Our result is the following:

\proclaim{Theorem 8}
If $b \in \BMO$, $1 < p < \infty$ and the inequalities (7), (8), (9) hold
(with $\gamma_i \geq 0$, $i = 1, \dots, N$) then there exists some constant $C
> 0$ such that
$$
   \Vert u [M_b, S_n](v^{-1} f) \Vert_{L^p(d\nu)}
   \leq C \Vert f \Vert_{L^p(d\nu)}
$$
for each $n \geq 0$ and $f \in L^p(d\nu)$.
\endproclaim

For the proof of Theorem 8 we firstly state the following result:

\proclaim{Lemma 9}
With the hypothesis of Theorem 8, we have:
\roster
\item"a)" $u(a_i) \Vert v^{-1} L_n(x,a_i) [b(x) - b(a_i)]
\Vert_{L^{p'}(d\mu)} \leq C, \qquad n \geq 0,\ i = 1, \dots, k$.

\item"b)" $v(a_i)^{-1} \Vert u L_n(x,a_i) [b(x) - b(a_i)] \Vert_{L^p(d\mu)}
\leq C, \qquad n \geq 0,\ i = 1, \dots, k$.
\endroster
\endproclaim

\demo{Proof}
Take $r$, $s$ such that $1/r + 1/s = 1/p'$. By H\"older's
inequality,
$$
   \Vert v^{-1} L_n(x,a_i) [b(x) - b(a_i)] \Vert_{L^{p'}(d\mu)} \leq
   \Vert v^{-1} L_n(x,a_i) \Vert_{L^{r}(d\mu)}
   \Vert [b(x) - b(a_i)] \Vert_{L^{s}(d\mu)}.
$$
From the John-Nirenberg inequality (13) it follows
$$
   \Vert b(x) - b(a_i) \Vert_{L^{s}(d\mu)} \leq C,
$$
while (10), (11), (12) lead to
$$
   \Vert v^{-1} L_n(x,a_i) \Vert_{L^{r}(d\mu)} \leq C
$$
provided $r$ is near enough to $p'$. This proves $a)$. Part $b)$ follows in a
similar way. \qed
\enddemo

Now, according to Lemma 3, we only need to prove
$$
   \Vert u [M_b, T_n](v^{-1}f) \Vert_{L^p(d\mu)}
   \leq C \Vert f \Vert_{L^p(d\mu)},
   \qquad f \in L^p(d\mu), \ n \geq 0.
$$
From the analog of Lemma 4, it is enough to show
$$
   \Vert u w^A [M_b,{\widetilde S}_n^A]([v w^A]^{-1} f) \Vert _{L^p(d\mu^A)}
   \leq C \Vert f \Vert _{L^p(d\mu^A)}
$$
for each $A \subseteq \{a_1,\dots,a_k\}$. For the sake of simplicity, we will
prove this inequality only for $A = \emptyset$. Then, $w^A = 1$, $d\mu^A =
d\mu$. Let us denote ${\widetilde S}_n = {\widetilde S}_n^\emptyset$.

\proclaim{Lemma 10}
With the hypothesis of Theorem 8, we have
$$
   \Vert u [M_b,{\widetilde S}_n](v^{-1} f) \Vert _{L^p(d\mu)}
   \leq C \Vert f \Vert _{L^p(d\mu)}
$$
for each $n \geq 0$, $f \in L^p(d\mu)$.
\endproclaim

\demo{Proof}
Let us denote by $K_n(x,y)$ the $n$-th kernel relative to $d\mu =
w(x) dx$, by ${\widetilde P}_n$ the orthonormal polynomials relative to $d\mu$
and by ${\widetilde Q}_n$ the orthonormal polynomials relative to
$(1-x^2)d\mu$. Then,
$$
   {\widetilde S}_n g(x) = \int_I K_n(x,y) g(y) w(y) dy.
$$
Also,
$$
   |{\widetilde P}_n(x) |
   \leq C (1 - x + n^{-2})^{-({\alpha \over 2} + {1 \over 4})}
   (1 + x + n^{-2})^{-({\beta \over 2} + {1 \over 4})}
   \prod_{i=1}^N (|x - t_i| + n^{-1})^{- {\gamma_i \over 2}},
$$
$$
   |{\widetilde Q}_n(x) |
   \leq C (1 - x + n^{-2})^{-({\alpha \over 2} + {3 \over 4})}
   (1 + x + n^{-2})^{-({\beta \over 2} + {3 \over 4})}
   \prod_{i=1}^N (|x - t_i| + n^{-1})^{- {\gamma_i \over 2}}
$$
(see \cite{B}). Now, we have Pollard's decomposition of $K_n$ (see \cite{P},
\cite{M 1}):
$$
   \eqalign{
   K_n(x,y) &= r_n {\widetilde P}_{n+1}(x) {\widetilde P}_{n+1}(y)
   \cr
   &+ s_n {\widetilde P}_{n+1}(x) {(1-y^2) {\widetilde Q}_n(y) \over x - y}
   \cr
   &+ s_n (1-x^2) {\widetilde Q}_n(x) {{\widetilde P}_{n+1}(y) \over y - x},
   \cr}
$$
for some bounded sequences $\{r_n\}$, $\{s_n\}$ of real numbers. Actually,
from $\mu' > 0$ a.e., it follows $\lim r_n = -1/2$, $\lim s_n = 1/2$
(it can be deduced from \cite{P} and either \cite{R} or \cite{MNT}). Therefore,
we can write
$$
   [M_b, {\widetilde S}_n]
   = r_n \Psi_{1,n} - r_n \Psi_{2,n} + s_n \Psi_{3,n} - s_n \Psi_{4,n},
$$
where
$$
   \Psi_{1,n} g(x) =
   [b(x) - b_I] {\widetilde P}_{n+1}(x)
   \int_I {\widetilde P}_{n+1}(y) g(y) w(y) dy,
$$
$$
   \Psi_{2,n} g(x) =
   {\widetilde P}_{n+1}(x)
   \int_I [b(y) - b_I] {\widetilde P}_{n+1}(y) g(y) w(y) dy,
$$
$$
   \Psi_{3,n} g(x) =
   {\widetilde P}_{n+1}(x)\ [M_b, H]((1-y^2) {\widetilde Q}_n g w)(x),
$$
$$
   \Psi_{4,n} g(x) =
   (1-x^2) {\widetilde Q}_n(x)\ [M_b, H]({\widetilde P}_{n+1} g w)(x).
$$
Lemma 11 below shows that the operators $\Psi_{3,n}$ are uniformly bounded; the
proof for $\Psi_{4,n}$ is entirely similar. Lemma 12 shows that the operators
$\Psi_{1,n}$ are uniformly bounded and the proof for $\Psi_{2,n}$ is again
similar. \qed
\enddemo

\proclaim{Lemma 11}
With the hypothesis of Theorem 8, there exists a constant $C > 0$ such that
$$
   \Vert u \Psi_{3,n}(v^{-1} f) \Vert_{L^p(d\mu)}
   \leq C \Vert f \Vert_{L^p(d\mu)}
$$
for each $n \geq 0$ and $f \in L^p(d\mu)$.
\endproclaim

\demo{Proof}
According to the definition of $\Psi_{3,n}$, we must show
$$
   \Vert [M_b, H] g \Vert_{L^p(u^p |{\widetilde P}_{n+1}|^p w)} \leq C
   \Vert g \Vert_{L^p(v^p |{\widetilde Q}_n|^{-p} (1-x^2)^{-p} w^{1-p})}.
$$
By the result of S. Bloom (\cite{Bl}), it is enough to find weights
$\{\phi_n\}$ and positive constants $K_1, K_2 > 0$ such that:
\roster
\item"a)" $K_1 u(x)^p |{\widetilde P}_{n+1}(x)|^p w(x) \leq \phi_n(x) \leq K_2
v(x)^p |{\widetilde Q}_n(x)|^{-p} (1-x^2)^{-p} w(x)^{1-p}$;

\item"b)" $\phi_n \in A_p(-1,1)$, with an $A_p$-constant independent of $n$.
\endroster

We take $\phi_n$ of the form
$$
   \eqalign{
   \phi_n(x) &= (1 - x)^{r_0} (1 - x + n^{-2})^{s_0}
   \cr
   &\times \prod_{i=1}^N |x - t_i|^{r_i} (|x - t_i| + n^{-1})^{s_i}
   \cr
   &\times (1 + x)^{r_{N+1}} (1 + x + n^{-2})^{s_{N+1}}
   \cr}
$$
Then, condition $b)$ is equivalent to
$$
   -1 < r_i < p - 1, \qquad -1 < r_i + s_i < p - 1
   \qquad (i = 0, 1, \dots, N+1)
   \tag 14
$$
(see \cite{GPV 2}). Now, it is not difficult to see from (7), (8), (9) and
$\gamma_i \geq 0$ that we can take $r_i$ such that
$$
   -1 < r_i < p-1 \qquad (i = 0, 1, \dots, N+1)
$$
$$
   A p - p + \alpha (1 - p) \leq r_0 \leq a p + \alpha
$$
$$
   B p - p + \beta (1 - p) \leq r_{N+1} \leq b p + \beta
$$
$$
   G_i p + \gamma_i (1 - p) \leq r_i \leq g_i p + \gamma_i
   \qquad (i = 1, \dots, N)
$$
and then $s_i$ such that
$$
   -1 < r_i + s_i < p-1 \qquad (i = 0, 1, \dots, N+1)
$$
$$
   A p - p + \alpha (1 - p) + p \left({\alpha \over 2} + {3 \over 4}\right)
   \leq r_0 + s_0 \leq
   a p + \alpha - p \left({\alpha \over 2} + {1 \over 4}\right)
$$
$$
   B p - p + \beta (1 - p) + p \left({\beta \over 2} + {3 \over 4}\right)
   \leq r_{N+1} + s_{N+1} \leq
   b p + \beta - p \left({\beta \over 2} + {1 \over 4}\right)
$$
$$
   G_i p + \gamma_i (1 - p) + p {\gamma_i \over 2}
   \leq r_i + s_i \leq
   g_i p + \gamma_i - p {\gamma_i \over 2} \qquad (i = 1, \dots, N).
$$
Using the estimates for ${\widetilde P}_n$ and ${\widetilde Q}_n$ and (14), we
can see that these conditions imply $a)$ and $b)$. \qed
\enddemo

\proclaim{Lemma 12}
With the hypothesis of Theorem 8, there exists a constant $C > 0$ such that
$$
   \Vert u \Psi_{1,n}(v^{-1} f) \Vert_{L^p(d\mu)}
   \leq C \Vert f \Vert_{L^p(d\mu)}
$$
for each $n \geq 0$ and $f \in L^p(d\mu)$.
\endproclaim

\demo{Proof}
Taking any $r > p$ and applying H\"older's inequality twice and
then (13), it follows
$$
   \Vert u \Psi_{1,n}(v^{-1} f) \Vert_{L^p(d\mu)}
   \leq C \Vert b \Vert_{\BMO}
   \ \Vert u {\widetilde P}_{n+1} w^{1/p} \Vert_{L^r(dx)}
   \ \Vert v^{-1} {\widetilde P}_{n+1} \Vert_{L^{p'}(w)}
   \ \Vert f \Vert_{L^p(w)}.
$$
Therefore, it is enough to have, for some $r > p$,
$$
   \Vert u {\widetilde P}_{n+1} w^{1/p} \Vert_{L^r(dx)}
   \ \Vert v^{-1} {\widetilde P}_{n+1} \Vert_{L^{p'}(w)} < C.
$$
This can be verified using the estimates for ${\widetilde P}_n$. \qed
\enddemo

\head 4. Weak behaviour for Jacobi weights with mass points on the
interval $[-1,1]$.\endhead

Let us consider now, as a particular case, a measure of the form
$$
   d\nu = (1-x)^{\alpha}(1+x)^{\beta}dx + \sum_{i=1}^k M_i\delta_{a_i}
$$
and $u=v=1$. If either $\alpha>-1/2$ or $\beta>-1/2$, then Corollary 6
determines an open interval of mean convergence $(p_0,p_1)$, where
$1<p_0<p_1<\infty$, and the $S_n$ are not uniformly bounded in $L^{p_0}(d\nu)$
or $L^{p_1}(d\nu)$. By symmetry, we can suppose $\alpha>-1/2$, $\alpha \geq
\beta >-1$, so that
$$
   p_0={4(\alpha+1) \over 2\alpha+3}, \qquad p_1={4(\alpha+1) \over 2\alpha+1}.
$$
In the absolutely continuous case and $\alpha=\beta=0$ (i.e.\ for Legendre
polynomials), when the mean convergence interval is $(4/3,4)$, Chanillo proved
(see \cite{C}) that the partial sums $S_n$ are not of weak type for $p=4$, that
is, there exists no constant $C>0$ such that for every $n \geq 0$ and $f \in
L^4(dx)$
$$
   \Vert S_nf \Vert _{L^{4,\infty}(dx)} \leq C \Vert f \Vert _{L^4(dx)}.
$$
It was also shown that these operators are of restricted weak type for $p=4$
(and $p=4/3$, by duality), that is, the previous inequality is verified if we
replace the $L^4$ norm by the $L^{p,1}$ norm. Actually, this is equivalent to
the inequality
$$
   \Vert S_n \chi_E \Vert _{L^{4,\infty}(dx)}
   \leq C \Vert \chi_E \Vert_{L^4(dx)}
$$
for every measurable set $E$ (see \cite{SW, Theorem 3.13}). The authors
obtained (see \cite{GPV 1}, \cite{GPV 2}) similar results for Jacobi weights.
The weak boundedness at the end points has also been considered for other
operators of Fourier Analysis. An important previous paper on the subject is
due to Kenig and Tomas (\cite{KT}), who studied the disk multiplier for radial
functions.

We can now prove that these results also hold with mass points:

\proclaim{Theorem 13}
Let $\alpha>-1/2$, $\alpha \geq \beta >-1$. If $p={4(\alpha+1) \over
2\alpha+1}$ or $p={4(\alpha+1) \over 2\alpha+3}$, then there exists no
constant $C$ such that for every $n \geq 0$ and $f \in L^p(d\nu)$
$$
   \Vert S_nf \Vert _{L^{p,\infty}(d\nu)} \leq C \Vert f \Vert _{L^p(d\nu)}.
$$
\endproclaim

\proclaim{Theorem 14}
Under the hypothesis of Theorem 13, there exists a constant $C>0$ such that
for every measurable set $E\subseteq [-1,1]$ and every $n \geq 0$
$$
   \Vert S_n\chi_E \Vert _{L^{p,\infty}(d\nu)} \leq C \Vert \chi_E \Vert
   _{L^p(d\nu)}.
$$
\endproclaim

Let us take now
$$
   u(x)=(1-x)^a(1+x)^b, \quad x\in (-1,1);
$$
$$
   0 < u(\pm 1) < \infty.
$$
We can extend Theorems 13 and 14 to the weighted case, when both $\alpha$ and
$\beta$ are greater or equal to $-1/2$.

\proclaim{Theorem 15}
Let $\alpha, \beta \geq -1/2$, $1<p<\infty$. If there exists a constant $C>0$
such that
$$
   \Vert uS_n(u^{-1} f) \Vert _{L^{p,\infty}(d\nu)}
   \leq C \Vert f \Vert_{L^p(d\nu)}
$$
for every $n \geq 0$ and $f \in L^p(d\nu)$, then the inequalities
$$
   \left| a + (\alpha + 1) \left( {1 \over p} - {1 \over 2} \right) \right|
   < {1 \over 4},
   \qquad
   \left| b + (\beta + 1) \left( {1 \over p} - {1 \over 2} \right) \right|
   < {1 \over 4}
$$
are verified.
\endproclaim

\proclaim{Theorem 16}
Let $\alpha, \beta \geq -1/2$, $1<p<\infty$. If the inequalities
$$
   -{1 \over 4} \leq a+(\alpha+1)\left({1 \over p}-{1 \over 2}\right)
   < {1 \over 4},
   \qquad
   -{1 \over 4} \leq b+(\beta+1)\left({1 \over p}-{1 \over 2}\right)
   < {1 \over 4}
$$
hold, then there exists a constant $C>0$ such that
$$
   \Vert uS_n(u^{-1}\chi_E) \Vert _{L^{p,\infty}(d\nu)}
   \leq C \Vert \chi_E \Vert_{L^p(d\nu)}
$$
for every $n\geq 0$ and every measurable set $E\subseteq [-1,1]$.
\endproclaim

\remark{Remark}
By standard arguments of duality (see \cite{GPV 2}), Theorem 16
also holds when
$$
   -{1 \over 4} < a+(\alpha+1)\left({1 \over p}-{1 \over 2}\right)
   \leq {1 \over 4},
   \qquad
   -{1 \over 4} < b+(\beta+1)\left({1 \over p}-{1 \over 2}\right)
   \leq {1 \over 4}.
$$
\endremark

\demo{Proof of Theorems 14 and 16}
We only need to show that
conditions $a)$, $b)$ and $c)$ in Lemma 1 are verified, with $p = q$, $r =
\infty$, $s = 1$. The estimates (10), (11), (12) are now
$$
   |L_n(x,a_i)|
   \leq C (1-x+n^{-2})^{-(2\alpha+1)/4}(1+x+n^{-2})^{-(2\beta+1)/4}
$$
if $a_i \neq \pm 1$;
$$
   |L_n(x,1)| \leq C (1+x+n^{-2})^{-(2\beta+1)/4}
$$
if $1$ is a mass point; and
$$
   |L_n(x,-1)| \leq C (1-x+n^{-2})^{-(2\alpha+1)/4}
$$
if $-1$ is a mass point, with $C$ independent of $n\geq 0$ and $x\in [-1,1]$.
Since either
$$
   (1-x+n^{-2})^{-(2\alpha+1)/4} \leq C
$$
or
$$
   (1-x+n^{-2})^{-(2\alpha+1)/4} \leq (1-x)^{-(2\alpha+1)/4},
$$
conditions $b)$ and $c)$ can be checked out taking into account that
$$
   (1-x)^r \in L^{p,\infty}((1-x)^sdx) \iff
   pr+s+1 \geq 0, (r,s) \ne (0,-1)
$$
(see \cite{GPV 2}, for example).

In order to prove condition $a)$, we can use Pollard's decomposition for the
kernels $L_n$ (see \cite{P}, \cite{M 1}) and write
$$
   T_nf = r_n W_{1,n}f + s_n W_{2,n}f - s_n W_{3,n}f,
$$
with
$$
   \eqalign{
   &\lim r_n = -1/2, \quad \lim s_n = 1/2,
   \cr
   &W_{1,n}f(x) = P_{n+1}(x) \int_{-1}^1P_{n+1}(y)f(y)w(y)dy,
   \cr
   &W_{2,n}f(x) = P_{n+1}(x) H((1-y)^2Q_n(y)f(y)w(y),x),
   \cr
   &W_{3,n}f(x) = (1-x^2)Q_n(x) H(P_{n+1}(y)f(y)w(y),x),
   \cr}
$$
where $\{P_n\}_{n\geq 0}$ is the sequence of orthonormal polynomials relative
to $d\nu$, $\{Q_n\}_{n\geq 0}$ is the sequence associated to $(1-x^2)d\nu$,
$w(x)=(1-x)^{\alpha}(1+x)^{\beta}$ and $H$ is the Hilbert transform on the
interval $[-1,1]$.

The polynomials $\{P_n\}_{n\geq 0}$ and $\{Q_n\}_{n\geq 0}$ have the
estimates (see \cite{GPRV 3})
$$
   |P_n(x)|
   \leq C (1-x+n^{-2})^{-(2\alpha+1)/4}(1+x+n^{-2})^{-(2\beta+1)/4}
$$
and
$$
   |Q_n(x)|
   \leq C (1-x+n^{-2})^{-(2\alpha+3)/4}(1+x+n^{-2})^{-(2\beta+3)/4}
$$
with $C$ independent of $n\geq 0$ and $x\in [-1,1]$, as in the absolutely
continuous case. We can now proceed exactly as in this case and show that
$W_{1,n}$, $W_{2,n}$ and $W_{3,n}$ are of restricted weak type (see \cite{GPV
2}).
\qed
\enddemo

\demo{Proof of Theorems 13 and 15}
If
$$
   \Vert uS_n(u^{-1} f) \Vert _{L^{p,\infty}(d\nu)}
   \leq C \Vert f \Vert_{L^p(d\nu)}
$$
for every $n \geq 0$ and $f \in L^p(d\nu)$, then
$$
   u \in L^{p,\infty}(d\nu),
$$
$$
   u^{-1} \in L^{p'}(d\nu),
$$
$$
   u w^{-1/2}(1-x^2)^{-1/4} \in L^{p,\infty}(w)
$$
and
$$
   u^{-1}w^{-1/2}(1-x^2)^{-1/4} \in L^{p'}(w)
$$
with $w(x)=(1-x)^{\alpha}(1+x)^{\beta}$ (see \cite{GPV 1}). This
proves Theorem 13 for $p={4(\alpha+1) \over 2\alpha+3}$ and implies
$$
   -{1 \over 4} \leq a+(\alpha+1)\left({1 \over p}-{1 \over 2}\right)
   < {1 \over 4},
   \quad
   -{1 \over 4} \leq b+(\beta+1)\left({1 \over p}-{1 \over 2}\right)
   < {1 \over 4}
$$
in Theorem 15. For
$$
   -{1 \over 4} = a+(\alpha+1)\left({1 \over p}-{1 \over 2}\right)
$$
or
$$
   -{1 \over 4} = b+(\beta+1)\left({1 \over p}-{1 \over 2}\right)
$$
in Theorem 15 and $p={4(\alpha+1) \over 2\alpha+1}$ in Theorem 13, it can be
proved that $W_{1,n}$ and $W_{3,n}$ are of weak type, while $W_{2,n}$ is not
(like in \cite{GPV 1}, \cite{GPV 2}). \qed
\enddemo

\head 5. Laguerre weights with a positive mass on $0$.\endhead

Lemma 1 is also useful to study the mean boundedness of the
Fourier series in the polynomials orthonormal with respect to the measure
$$
   d\nu = e^{-x}x^{\alpha}dx + M\delta_0
$$
on $[0,\infty)$ (that is, a Laguerre weight with a mass $M>0$ at $0$). In this
case, parts $b)$ and $c)$ can be handled having in mind that the kernels
$L_n(x,0)$ admit the formula
$$
   L_n(x,0) = r_n Q_n(x),
   \tag 15
$$
where $Q_n$ is the $n$-th orthonormal Laguerre polynomial relative to the
measure $e^{-x}x^{\alpha+1}dx$. This formula follows from the fact that
$$
   \int_0^{\infty} L_n(x,0) xR_{n-1}(x) [e^{-x}x^{\alpha}dx + M\delta_0] = 0
$$
and
$$
   \int_0^{\infty} Q_n(x) xR_{n-1}(x) [e^{-x}x^{\alpha}dx + M\delta_0] = 0
$$
for any polynomial $R_{n-1}$ of degree at most $n-1$. The constants
$r_n = L_n(0,0)/Q_n(0)$ can be asymptotically estimated. As we mentioned in the
proof of Lemma 1, $\{L_n(0,0)\}$ is an increasing, bounded sequence, since $0$
is a mass point (see \cite{N, p. 4}). On the other hand, if we denote by
$\{L_n^{\alpha+1}\}$ the classical, not normalized Laguerre polynomials
relative to $e^{-x}x^{\alpha+1}dx$, then it is well known that 
$$
   L_n^{\alpha+1}(0) = {\Gamma(n+\alpha+2) \over n!\Gamma(\alpha+2)}
$$
(see \cite{S}, \cite{M 3}) and
$$
   \Vert L_n^{\alpha+1} \Vert _{L^2(e^{-x}x^{\alpha+1}dx)}
   = {\Gamma(n+\alpha+2) \over n!},
$$
what, with our notation, implies
$$
   Q_n(0) = {\Gamma(n+\alpha+2)^{1/2} \over \Gamma(\alpha+2) (n!)^{1/2}}
   \sim n^{(\alpha+1)/2}.
$$
Therefore,
$$
   r_n \sim n^{-(\alpha+1)/2}.
   \tag 16
$$
According to (15) and (16), in order to find bounds for the kernels $L_n(x,0)$
we only need bounds for the normalized classical Laguerre polynomials. These
bounds, as well as boundedness results for Laguerre series, can be found in
Muckenhoupt's paper \cite{M 3}. Thus, we can use Lemma 1 as in the generalized
Jacobi case to find that Muckenhoupt's result (\cite{M 3, Theorem 7}) remains
valid in the case of a Laguerre weight with a positive mass on $0$. The same
can be done for Hermite series (see \cite{M 3, Theorem 1}).

%-----------------------------------------

\enddocument